# Correction to Connecting invariant manifolds and the solution of the $C^1$ stability and Ω-stability conjectures for flows

By Shuhei Hayashi

There is a gap in the proof of Lemma VII.4 in [1]. We present an alternative proof of Theorem B ($C^1$ Ω-stable vector fields satisfy Axiom A) in [1]. The novel and essential part in the proof of the stability and Ω-stability conjectures for flows is the connecting lemma introduced in [1]. A mistake in the proof of the last conjecture was pointed out to me by Toyoshiba [5], who later also provided an independent proof of it, again based on the connecting lemma and previous arguments by Mañé and Palis.

The crucial step to the proof of Theorem B is the separation of singularities from periodic orbits ([1, Corollary III]) by the $C^1$ connecting lemma ([1, Theorem A]). After the separation, the proof proceeds based on Mañé's theorems used in [3] and we still rely on Palis's argument in [4], proving first the density of Axiom A diffeomorphisms in the set of $C^1$ Ω-stable ones to then show that every $C^1$ Ω-stable diffeomorphism satisfies Axiom A.

Let $\mathcal{G}^1_\Omega(M)$ be the set of $C^1$ Ω-stable vector fields on a compact smooth boundaryless manifold $M$ with the $C^1$ topology and $X \in \mathcal{G}^1_\Omega(M)$. As in [1], we prove the hyperbolicity of $\overline{\mathrm{Per}(X)}$ ($= \Omega(X) - \mathrm{Sing}(X)$) by induction. In fact, we prove that $\overline{P}_j(X)$ is hyperbolic assuming that $\bigcup_{i=0}^{j-1} \overline{P}_i(X)$ is hyperbolic for some $1 \leq j \leq \dim M - 1$, where $\overline{P}_i(X)$ is the closure of the set of periodic points with index $i$ (dimension of the stable subspace), which is enough to conclude that $X$ satisfies Axiom A. For a dense subset of $\mathcal{G}^1_\Omega(M)$, we can use the statement of [1, Lemma VII.4] by an already classic argument on set-valued functions of $C^1$ vector fields. In fact, there is a residual subset of the set of $C^1$ vector fields (therefore of $\mathcal{G}^1_\Omega(M)$) in which the closure of the set of hyperbolic periodic points of saddle type moves continuously with respect to vector fields (see for instance the proof of [1, Corollary II] for this kind of argument). Therefore, as proved in [1], we get the density of Axiom A vector fields in $\mathcal{G}^1_\Omega(M)$. Then, by Ω-conjugacy, we see that $\Omega(X)$ can be decomposed into a finite union of disjoint compact invariant sets which are isolated and transitive. Moreover, Palis's argument ([4]) for flows shows that each component is homogeneous in the sense that the index of every periodic point in it is the same. Thus, the proof of Theorem B is reduced to proving the following claim:

*Claim.* Every homogeneous component of $\overline{P}_j(X)$ is hyperbolic.

Let $\mathcal{G}^1(M)$ be the interior of the set of $C^1$ vector fields on $M$, with the $C^1$ topology, such that all periodic orbits and singularities are hyperbolic. Then $\mathcal{G}^1_\Omega(M) \subset \mathcal{G}^1(M)$. Denote by $L^X_t$, $t \in \mathbf{R}$ the linear Poincaré flow of $X \in \mathcal{G}^1(M)$ on $N^*$ (see [1, p. 126] for the definition). As in [1, p. 131], let $N^*|\overline{P}_j(X) = E_j \oplus F_j$ be the dominated splitting such that

$$\|L^X_m|E_j(y)\| \cdot \|L^X_{-m}|F_j(X_m(y))\| \leq \lambda \tag{1}$$

for all $y \in \overline{P}_j(X)$ with $m \in \mathbf{Z}^+$ and $0 < \lambda < 1$ given by [1, Lemma VII.1], which is the continuous extension of hyperbolic splittings of periodic orbits of index $j$ with respect to $L^X_t$. To prove the Claim, it is enough to show that $E_j$ is contracting by the following lemma proved in [3, Theorem II.1], which is Lemma VII.5 in [1] and written for this setting:

LEMMA 1 (Mañé). *Let $\Lambda$ be a compact invariant set of $X \in \mathcal{G}^1(M)$ such that $\Lambda \cap \mathrm{Sing}(X) = \emptyset$ and $\Omega(X|\Lambda) = \Lambda$. Suppose that $N^*|\Lambda = E \oplus F$ is a dominated splitting such that the dimension of the subspaces $E(y)$, $y \in \Lambda$ is constant,*

$$\|L^X_m|E(y)\| \cdot \|L^X_{-m}|F(X_m(y))\| \leq \lambda$$

*for all $y \in \Lambda$, and*

$$\liminf_{n \to +\infty} \frac{1}{n} \sum_{\ell=1}^n \log \|L^X_{-m}|F(X_{m\ell}(x))\| \leq \log \lambda$$

*holds for a dense set of points $x \in \Lambda$, where $m \in \mathbf{Z}^+$ and $0 < \lambda < 1$ are given in (1). Then, if $E$ is contracting, $F$ is expanding (and therefore $\Lambda$ is hyperbolic).*

Let $\sum(X)$ be the set of "strongly closable points" given in [1, Lemma VII.6 (Ergodic Closing Lemma for time-one maps)] and originally introduced by Mañé in [2]. We shall need the following lemma:

LEMMA 2. *Let $X \in \mathcal{G}^1(M)$. If $x \in \sum(X)$ and $\overline{\mathcal{O}_X(x)} \cap \mathrm{Sing}(X) = \emptyset$, then $\overline{\mathcal{O}_X(x)}$ contains a hyperbolic set, where $\mathcal{O}_X(x) = \{X_t(x) : t \in \mathbf{R}\}$.*

*Proof.* We can suppose that $x \in \sum(X) - \mathrm{Per}(X)$. Let $\mathcal{U}_n$, $n \geq 0$, be a basis of neighborhoods of $X$. Then, by the definition of $\sum(X)$ ([1, p. 132]; see also [2, p. 506]), there exists $\{t_n > 0 : n \geq 0\}$ with $\lim_{n \to +\infty} t_n = +\infty$, $X^n \in \mathcal{U}_n$ and $y_n \in \mathrm{Per}(X^n)$ having period $T_n$ such that $\{X_t(x) : 0 \leq t \leq t_n\}$ can be approximated by $\{X^n_t(y_n) : 0 \leq t \leq T_n\}$ for large $n$. Without loss of generality we may assume that the index of the $X^n$-periodic point $y_n$ is the same for all $n \geq 0$ (by taking a subsequence if necessary). Then, by [1, Lemma VII.1], the following properties hold for all $n \geq 0$ with $T_n \geq m$:

$$\|L_m^{X^n}|E_n^s(y)\| \cdot \|L_{-m}^{X^n}|E_n^u(X_m^n(y))\| \leq \lambda$$

for all $y \in \{X_t^n(y_n) : 0 \leq t \leq T_n\}$, and

$$\prod_{\ell=1}^{[T_n/m]} \|L_{-m}^{X^n}|E_n^u(X_{m\ell}^n(y_n))\| \leq K\lambda^{[T_n/m]},$$

where $E_n^s \oplus E_n^u$ is the hyperbolic splitting with respect to $L_t^{X^n}$ over the $X^n$-periodic orbit with $y_n$. Note that the dimension of $E_n(y)$ is constant, and the angle between $E_n^s(y)$ and $E_n^u(y)$ is uniformly bounded away from 0 by [2, Lemma II.9]. Then, defining the splitting $E \oplus F$ over $\mathcal{O}_X(x)$ by accumulation of $\{E_n^s \oplus E_n^u : n \geq 0\}$ (see the definition of $\sum(X)$ again), we have, by continuity, the following properties:

$$\|L_m^X|E(y)\| \cdot \|L_{-m}^X|F(X_m(y))\| \leq \lambda$$

for all $y \in \mathcal{O}_X(x)$, and

$$\liminf_{n \to +\infty} \frac{1}{n} \sum_{\ell=1}^{n} \log \|L_{-m}^X|F(X_{m\ell}(x))\| \leq \log \lambda.$$

It is well known that the dominated splitting $E \oplus F$ over $\mathcal{O}_X(x)$ can be continuously extended to $\widetilde{E} \oplus \widetilde{F}$ over $\overline{\mathcal{O}_X(x)}$ with the same $m \in \mathbf{Z}^+$ and $0 < \lambda < 1$. Hence, we can apply Lemma 1 to $\Lambda = \overline{\mathcal{O}_X(x)}$. If $\overline{\mathcal{O}_X(x)}$ is not hyperbolic; that is, $\widetilde{E}$ is not contracting, then, as in [1, p. 132], there exists $p \in \overline{\mathcal{O}_X(x)} \cap \sum(X)$ such that

(2) $$\lim_{n \to +\infty} \frac{1}{n} \sum_{\ell=0}^{n-1} \log \|L_m^X|\widetilde{E}(X_{m\ell}(p))\| \geq 0.$$

When $p \in \text{Per}(X)$, $\mathcal{O}_X(p)$ is a hyperbolic set, we may assume that $p \notin \text{Per}(X)$. Then, we can continue this argument for $\overline{\mathcal{O}_X(p)}$ instead of $\overline{\mathcal{O}_X(x)}$. As observed in [1, pp. 132–133], the index $i_0$ of periodic point created by the Ergodic Closing Lemma from $\mathcal{O}_X(p)$ is less than $\dim \widetilde{E}$. If $i_0 = 0$, then, by the same argument as in the proof of the finiteness of periodic orbits in $\overline{P_0}(M)$, $p$ cannot be reccurent. Therefore, we may suppose that $\overline{\mathcal{O}_X(p)}$ has a contracting subbundle (by continuing this further if necessary) and therefore, by Lemma 1, $\overline{\mathcal{O}_X(p)}$ is hyperbolic, proving Lemma 2.

Now let us prove the Claim. Assume that a homogeneous component $\widetilde{\Lambda}$ of $\overline{P}_j(X)$ is not hyperbolic. Then, we can find $p \in \sum(X) \cap \widetilde{\Lambda}$ satisfying property

(2) with $\widetilde{E}$ replaced by $E_j$ in (1). Lemma 2 implies that $\overline{\mathcal{O}_X(p)}$ contains a hyperbolic set $\Lambda$, and, as in the proof of Lemma 2, the index of any periodic point in $\Lambda$ is less than $j$. This contradicts the homogeneity of $\widetilde{\Lambda}$.

SCHOOL OF COMMERCE, WASEDA UNIVERSITY, SHINJUKU, TOKYO, JAPAN
*E-mail address*: shuhei@mn.waseda.ac.jp

## REFERENCES

[1] S. HAYASHI, Connecting invariant manifolds and the solution of the $C^1$ stability and $\Omega$-stability conjectures for flows, Ann. of Math. **145** (1997), 81–137.
[2] R. MAÑÉ, An ergodic closing lemma, Ann. of Math. **116** (1982), 503–540.
[3] ______, A proof of the $C^1$ stability conjecture, Publ. Math. I.H.E.S. **66** (1988), 161–210.
[4] J. PALIS, On the $C^1$ $\Omega$-stability conjecture, Publ. Math. I.H.E.S. **66** (1988), 211–215.
[5] H. TOYOSHIBA, private communication.